\documentclass[a4paper,11pt]{article}
\usepackage[latin1]{inputenc}
\usepackage[T1]{fontenc}
\usepackage[francais]{babel}
\usepackage{amssymb}
\usepackage{amsfonts, amsmath}
\newtheorem{thm}{Th\'{e}or\`{e}me}

\title{Fonctions m\'{e}romorphes et fonctions th\^{e}ta sur les surfaces de Riemann}
\author{A. Lesfari
\\\emph{Department of Mathematics}
\\\emph{Faculty of Sciences}
\\\emph{University of Choua\"{i}b Doukkali}
\\\emph{B.P. 20, 24000 El Jadida, Morocco}.
\\\emph{lesfariahmed@yahoo.fr}}
\date{}

\begin{document}
\maketitle Theta functions play a major role in many current
researches and are powerful tools for studying integrable systems.
The purpose of this paper is to provide a short and quick
exposition of some aspects of meromorphic theta functions for
compact Riemann surfaces. The study of theta functions will be
done via an analytical approach using meromorphic functions in the
framework of Mumford. Some interesting examples will be given :
the classical Kirchhoff equations in the cases of Clebsch and
Lyapunov-Steklov, the Landau-Lifshitz equation and the sine-Gordon
equation.\\
\emph{AMS 2010 Subject Classification}: 30F10, 30D30, 14K25.\\
\emph{Key words}: Riemann surfaces, meromorphic functions, theta
functions.

\section{Fonctions th\^{e}ta}

Soient $X$ une surface de Riemann compacte de genre $g$ et
$B=(b_{jk})_{1\leq j,k\leq g}$ une matrice carr\'{e}e d'ordre $g$,
sym\'{e}trique et dont la partie imaginaire est d\'{e}finie
positive. On consid\`{e}re la fonction th\^{e}ta $\theta(z|B)$ de
Riemann d\'{e}finie \`{a} l'aide de la s\'{e}rie :
\begin{equation}\label{eqn:euler}
\theta(z|B)=\sum_{m\in \mathbb{Z}^g}e^{\pi i\langle
Bm,m\rangle+2\pi i\langle z,m\rangle},\quad z\in\mathbb{C}^g
\end{equation}
o\`{u} $\displaystyle{\langle
Bm,m\rangle=\sum_{j,k=1}^gb_{jk}m_jm_k}$, $\displaystyle{\langle
z,m\rangle=\sum_{j=1}^gz_jm_j}$. La convergence de cette s\'{e}rie
pour tout $z\in\mathbb{C}^g$, r\'{e}sulte du fait que $\mbox{Im
}B\geq 0$. On montre que cette s\'{e}rie converge absolument et
uniform\'{e}ment sur des ensembles compacts et qu'en outre, la
fonction $\theta(z|B)$ est holomorphe sur $\mathbb{C}^g$. On
posera dans la suite $\theta(z)\equiv\theta(z|B)$ lorsque la
matrice $B$ est fix\'{e}e. Soit $(e_1,...,e_g)$ une base de
$\mathbb{C}^g$ avec $(e_j)_k=\delta_{jk}$, et d\'{e}signons par $f_j=(b_{1j} ... b_{gj})^\top$, les colonnes de la matrice $B$ ou sous
forme condens\'{e}e $f_j=Be_j$, $j=1,...,g$.
\begin{thm}
La fonction $\theta$ satisfait aux \'{e}quations fonctionnelles :
\begin{equation}\label{eqn:euler}
\theta(z+e_j)=\theta(z),\qquad
\theta(z+f_j)=e^{-\pi ib_{jj}-2\pi iz_j}.\theta(z).
\end{equation}
Pour tout $m,n\in \mathbb{Z}^g$, on a
\begin{equation}\label{eqn:euler}
\theta(z+n+Bm)=e^{-\pi i\langle Bm,m\rangle-2\pi i\langle
m,z\rangle}.\theta(z).
\end{equation}
Les vecteurs de la forme $n+Bm$ forment un r\'{e}seau de
p\'{e}riodes.
\end{thm}
\emph{D\'{e}monstration}: La premi\`{e}re relation r\'{e}sulte de
la formule (1). Concernant la seconde relation, on a
\begin{eqnarray}
\theta(z+f_j&=&\sum_{m\in \mathbb{Z}^g}e^{\pi i\langle
Bm,m\rangle+2\pi i\langle m,z+f_j\rangle},\nonumber\\
&=&\sum_{m\in \mathbb{Z}^g}e^{\pi i\langle
B(n-e_j),n-e_j\rangle+2\pi i\langle n-e_j,z+f_j\rangle},\quad n\equiv m+e_j\nonumber\\
&=&e^{-\pi i\langle Be_j,e_j\rangle-2\pi
i\langle e_j,z\rangle}.\theta(z),\nonumber\\
&=&e^{-\pi ib_{jj}-2\pi iz_j}.\theta(z),\nonumber
\end{eqnarray}
et la, relation (3) en r\'{e}sulte imm\'{e}diatement. $\square$

Les vecteurs $e_1,...,e_g$ forment une base de p\'{e}riodes de la
fonction $\theta(z)$. Les vecteurs $f_j$ s'appellent les
quasi-p\'{e}riodes de $\theta(z)$. La fonction $\theta$ est
quasi-p\'{e}riodique et elle est bien d\'{e}finie sur la
vari\'{e}t\'{e} jacobienne de $X$. Consid\'{e}rons maintenant une
g\'{e}n\'{e}ralisation de la fonction th\^{e}ta (1) appelée
fonction th\^{e}ta de caract\'{e}ristiques
$\left[\begin{array}{c}\alpha\\\beta\end{array}\right]$
d\'{e}finie par
\begin{eqnarray}
\theta\left[\begin{array}{c}\alpha\\\beta\end{array}\right](z|B)&=&\sum_{m\in
\mathbb{Z}^g}e^{\pi i\langle B(m+\alpha),m+\alpha\rangle+2\pi
i\langle z+\beta,m+\alpha\rangle}, \quad \alpha,\beta\in \mathbb{R}^g\\
&=&e^{\pi i\langle B\alpha,\alpha\rangle+2\pi i\langle
z+\beta,\alpha\rangle}.\theta(z+\beta+B\alpha).
\end{eqnarray}
Pour all\'{e}ger les formules, on notera simplement :
$\theta\left[\begin{array}{c}\alpha\\\beta\end{array}\right](z)
\equiv\theta\left[\begin{array}{c}\alpha\\\beta\end{array}\right](z|B)$
lorsque la matrice $B$ est fix\'{e}e. En particulier,
$\theta\left[\begin{array}{c}0\\0\end{array}\right](z)=\theta(z)$.
D'apr\`{e}s la relation (3), on a aussi
$\theta\left[\begin{array}{c}m\\n\end{array}\right](z)=\theta(z)$,
$m,n\in \mathbb{Z}^g$. D\`{e}s lors, il suffit de consid\'{e}rer
les fonctions
$\theta\left[\begin{array}{c}\alpha\\\beta\end{array}\right](z)$
o\`{u} $\alpha=(\alpha_1,...,\alpha_g)$,
$\beta=(\beta_1,...,\beta_g)\in\mathbb{R}^g$ sont tels que :
$0<\alpha_j,\beta_j<1$, $j=1,...,g$.
\begin{thm}
La fonction $\theta$ v\'{e}rifie la propri\'{e}t\'{e} de
p\'{e}riodicit\'{e} suivante :
$$\theta\left[\begin{array}{c}\alpha\\\beta\end{array}\right](z)(z+n+Bm)=
e^{-\pi i\langle Bm,m\rangle-2\pi i\langle z,m\rangle+2\pi
i(\langle \alpha,n\rangle-\langle
\beta,m\rangle)}.\theta\left[\begin{array}{c}\alpha\\\beta\end{array}\right](z).$$
\end{thm}
\emph{D\'{e}monstration}: Il suffit d'utiliser un raisonnement
similaire \`{a} celui de la proposition pr\'{e}c\'{e}dente.
$\square$.

Lorsque $\alpha_1,...,\alpha_g$ et $\beta_1,...,\beta_g$ ne
prennent que les valeurs $0$ ou $\frac{1}{2}$, on dira que
$\left[\begin{array}{c}\alpha\\\beta\end{array}\right]$ est une
demi-p\'{e}riode. En outre, une demi-p\'{e}riode
$\left[\begin{array}{c}\alpha\\\beta\end{array}\right]$ est dite
paire si $4\langle\alpha,\beta\rangle\equiv 0$ (mod. 2) et impaire
sinon.

\begin{thm}
La fonction
$\theta\left[\begin{array}{c}\alpha\\\beta\end{array}\right](z)$
est paire (resp. impaire) si la demi-p\'{e}riode
$\left[\begin{array}{c}\alpha\\\beta\end{array}\right]$ est paire
(resp. impaire). En outre, on a $\theta(z)=\theta(-z)$.
\end{thm}
\emph{D\'{e}monstration}: En faisant la substitution $z\longmapsto
-z$, $m\longmapsto -m-2$, dans la relation (4), on obtient
imm\'{e}diatement pour le terme g\'{e}n\'{e}ral de la s\'{e}rie,
$$e^{\pi i\langle B(-m-\alpha),-m-\alpha\rangle+2\pi i\langle
-z+\beta,-m-\alpha\rangle}=e^{\pi i\langle
B(m+\alpha),m+\alpha\rangle+2\pi i\langle
z+\beta,m+\alpha\rangle}.e^{4\pi i\langle \alpha,\beta\rangle}.$$
Or d'apr\`{e}s la d\'{e}finition ci-dessus, le signe de $e^{4\pi
i\langle \alpha,\beta\rangle}$ est d\'{e}termin\'{e} par la
parit\'{e} du nombre $4\langle \alpha,\beta\rangle$, et la derni\`{e}re relation en r\'{e}sulte. $\square$

Par exemple le nombre de demi-p\'{e}riodes paires est \'{e}gal
\`{a} $2^{g-1}(2^g+1)$ et celui de demi-p\'{e}riodes impaires
\`{a} $2^{g-1}(2^g-1)$.

\section{Fonctions m\'{e}romorphes exprim\'{e}es en termes de fonctions th\^{e}ta}

Consid\'{e}rons le cas des surfaces de Riemann de genre $1$,
c.-\`{a}-d., des courbes elliptiques. Rappelons qu'une fonction
elliptique est une fonction m\'{e}romorphe doublement
p\'{e}riodique. La matrice $B$ se r\'{e}duit dans ce cas \`{a} un
nombre que l'on note $b$ avec $\mbox{Im }B\geq 0$. Les nombres $1$
et $b$ engendrent un parall\'{e}logramme des p\'{e}riodes not\'{e}
$\Omega$. Les quatre fonctions th\^{e}ta correspondant aux
demi-p\'{e}riodes
$\left[\begin{array}{c}1/2\\1/2\end{array}\right]$,
$\left[\begin{array}{c}1/2\\0\end{array}\right]$,
$\left[\begin{array}{c}0\\0\end{array}\right]$,
$\left[\begin{array}{c}0\\1/2\end{array}\right]$ sont
$i\theta_1(z)\equiv\theta\left[\begin{array}{c}1/2\\1/2\end{array}\right](z)$,
$\theta_2(z)\equiv\theta\left[\begin{array}{c}1/2\\0\end{array}\right](z)$,
$\theta_3(z)\equiv\theta\left[\begin{array}{c}0\\0\end{array}\right](z)=\theta(z)$,
$\theta_4(z)\equiv\theta\left[\begin{array}{c}0\\1/2\end{array}\right](z)$.
Ces fonctions sont holomorphes sur $\mathbb{C}$. En outre, on
d\'{e}duit imm\'{e}diatement du th\'{e}or\`{e}me 3 que
$\theta_1(z)$ est impaire et que $\theta_2(z)$, $\theta_3(z)$,
$\theta_4(z)$ sont paires. Pour d\'{e}terminer les z\'{e}ros des
fonctions $\theta_j$, il suffit d'apr\`{e}s le th\'{e}or\`{e}me 2
de les chercher dans le parall\'{e}logramme des p\'{e}riodes
$\Omega$. Comme $\theta_1(z)$ est impaire, alors $\theta_1(0)=0$
et les autres z\'{e}ros de $\theta_j(z)$ s'obtiennent via le
th\'{e}or\`{e}me 2. Prenons par exemple le cas de $\theta_3(z)$.
On a $\theta_3(\frac{1}{2}+\frac{b}{2})=0$ et $z=\frac{1}{2}(1+b)$
est le seul z\'{e}ro de cette fonction dans $\Omega$. En effet, on
a
\begin{eqnarray}
\frac{1}{2\pi
i}\int_{\partial\Omega}\frac{d\theta_3(z)}{\theta_3(z)}
&=&\frac{1}{2\pi i}\int_0^1(d\log\theta_3(z)-d\log\theta_3(z+b))\nonumber\\
&&\qquad+\frac{1}{2\pi
i}\int_0^b(d\log\theta_3(z+1)-d\log\theta_3(z)).\nonumber
\end{eqnarray}
D'apr\`{e}s le th\'{e}or\`{e}me 2, on a
$\theta_3(z+1)=\theta_3(z)$ et $\theta_3(z+b)=e^{-\pi ib-2\pi
iz}\theta_3(z)$, donc
$$\int_0^b(d\log\theta_3(z+1)-d\log\theta_3(z))=0,$$
et
\begin{eqnarray}
\int_0^1(d\log\theta_3(z)-d\log\theta_3(z+b))
&=&\int_0^1(d\log\theta_3(z)-d\log e^{-\pi ib-2\pi iz}\theta_3(z)),\nonumber\\
&=&\int_0^12\pi idz=1.\nonumber
\end{eqnarray}
D\`{e}s lors,
$$
\frac{1}{2\pi
i}\int_{\partial\Omega}\frac{d\theta_3(z)}{\theta_3(z)}
=\frac{1}{2\pi i}\int_{\partial\Omega}d\theta_3(z)=1.
$$
Par cons\'{e}quent, on a le r\'{e}sultat suivant :

\begin{thm}
La fonction $\theta(z)$ poss\`{e}de dans le parall\'{e}logramme
des p\'{e}riodes $\Omega$ (engendr\'{e} par $1$ et $b$), un seul
z\'{e}ro au point $z=\frac{1}{2}(1+b)$.
\end{thm}

En posant $z=x\in\mathbb{R}$, $b=it$, $t\in\mathbb{R}_+$,
l'\'{e}quation (1) s'\'{e}crit
$$
\theta(x|it)=\sum_{m\in\mathbb{Z}}e^{-\pi m^2t+2\pi
imx}=1+2\sum_{m=1}^\infty e^{-\pi m^2t}\cos 2\pi mx.
$$
Cette fonction est p\'{e}riodique par rapport \`{a} $x$,
c.-\`{a}-d., $\theta(x+1|it)=\theta(x|it)$, et satisfait \`{a}
l'\'{e}quation de la chaleur
$$4\pi \frac{\partial \theta(x|it)}{\partial t}=\frac{\partial^2 \theta(x|it)}{\partial
x^2}.$$ L'unicit\'{e} de cette solution r\'{e}sulte du fait que
$\displaystyle{\lim_{t\rightarrow
0}\theta(x|it)=\sum_{m=-\infty}^\infty \delta_m(x)}$, o\`{u}
$\delta_m$ est la distribution de Dirac au point $m$. De m\^{e}me,
la fonction $\theta_1(z)$ v\'{e}rifie une \'{e}quation
diff\'{e}rentielle de $3^{\mbox{\`{e}me}}$ ordre. En effet, il
suffit d'utiliser la relation
$\displaystyle{\wp(z)=-\frac{\partial^2}{\partial
z^2}\log\theta_1(z)+C}$, o\`{u} $C$ est une constante, $\wp(z)$
est la fonction de Weierstrass d\'{e}finie par
$$
\wp(z)=\frac{1}{z^2}+\sum_{\omega \in \Lambda \backslash
\{0\}}\left(\frac{1}{(z-\omega )^2} -\frac{1}{\omega ^2}\right),
$$
$\Lambda =\mathbb{Z}\omega _{1}\oplus\mathbb{Z}\omega _{2}$, est
le r\'{e}seau  engendr\'{e} par deux nombres complexes $\omega
_{1}$\ et $\omega _{2}$ diff\'{e}rents de $0$ tels que:
$\mbox{Im}\displaystyle{\left( \frac{\omega _{2}}{\omega
_{1}}\right)}>0$, et tenir compte de l'\'{e}quation
diff\'{e}rentielle :
\begin{equation}\label{eqn:euler}
\left(\wp ^{\prime
}(z)\right)^2=4\left(\wp(z)\right)^3-g_{2}\wp(z)-g_3,
\end{equation}
o\`{u} $g_{2}=60\displaystyle{\sum_{\omega \in \Lambda \backslash
\{0\}} \frac{1}{\omega ^{4}}}$,
$g_{3}=140\displaystyle{\sum_{\omega \in \Lambda \backslash
\{0\}}\frac{1}{\omega^{6}}}$. Par ailleurs, on a les identit\'{e}s classiques [11]:
\begin{thm}
La fonction th\^{e}ta satisfait aux formules d'addition :\\
$$\theta\left[\begin{array}{c}\alpha\\\beta\end{array}\right](z_1+z_2)
=\sum_{2\delta\in(\mathbb{Z}_2)^g}\widehat{\theta}\left[\begin{array}{c}
\frac{\alpha+\beta}{2}+\delta\\\\\gamma+\varepsilon\end{array}\right](2z_1).
\widehat{\theta}\left[\begin{array}{c}
\frac{\alpha-\beta}{2}+\delta\\\\\gamma-\varepsilon\end{array}\right](2z_2),$$
o\`{u} $\alpha,\beta,\gamma,\varepsilon\in\mathbb{R}^g$,
$$\theta\left[\begin{array}{c}\alpha\\\beta\end{array}\right](z)\equiv
\theta\left[\begin{array}{c}\alpha\\\beta\end{array}\right](z|B),\qquad
\widehat{\theta}\left[\begin{array}{c}\alpha\\\beta\end{array}\right](z)\equiv
\theta\left[\begin{array}{c}\alpha\\\beta\end{array}\right](z|2B).$$
$$
\theta\left[\begin{array}{c}m_1\\n_1\end{array}\right](z_1).
\theta\left[\begin{array}{c}m_2\\n_2\end{array}\right](z_2).
\theta\left[\begin{array}{c}m_3\\n_3\end{array}\right](z_3).
\theta\left[\begin{array}{c}m_4\\n_4\end{array}\right](z_4)
$$$$=\frac{1}{2^g}\sum_{2(a_1,a_2)\in(\mathbb{Z}_2)^{2g}}e^{-4\pi
i\langle
m_1,a_2\rangle}.\theta\left[\begin{array}{c}k_1+a_1\\l_1+a_2\end{array}\right](w_1)
...\theta\left[\begin{array}{c}k_4+a_1\\l_4+a_2\end{array}\right](w_4),$$
o\`{u} $(z_1,...,z_4)=(w_1,...,w_4)M$ avec
$${M}=\frac{1}{2}\left(\begin{array}{cccc}
\textbf{1}&\textbf{1}&\textbf{1}&\textbf{1}\\
\textbf{1}&\textbf{1}&-\textbf{1}&-\textbf{1}\\
\textbf{1}&-\textbf{1}&\textbf{1}&-\textbf{1}\\
\textbf{1}&-\textbf{1}&-\textbf{1}&\textbf{1}
\end{array}\right)
$$
Ici
$\left(\begin{array}{c}m_1\\n_1\end{array}\right)$,...,$\left(\begin{array}{c}m_4\\n_4\end{array}\right)$,
$\left(\begin{array}{c}k_1\\l_1\end{array}\right)$,...,$\left(\begin{array}{c}k_4\\l_4\end{array}\right)$
sont des vecteurs quelconques d'ordre $2g$ avec
$$\left(\left(\begin{array}{c}m_1\\n_1\end{array}\right),...,\left(\begin{array}{c}m_4\\n_4\end{array}\right)\right)
=
\left(\left(\begin{array}{c}k_1\\l_1\end{array}\right),...,\left(\begin{array}{c}k_4\\l_4\end{array}\right)\right)M,$$
et $\textbf{1}$ d\'{e}signe la matrice unit\'{e} d'ordre $g$ ou
$2g$.
\end{thm}
En particulier, on a les formules :
\begin{eqnarray}
\left(\theta\left[\begin{array}{c}0\\0\end{array}\right](z)\right)^2.
\left(\theta\left[\begin{array}{c}0\\0\end{array}\right](0)\right)^2
&=&
\left(\theta\left[\begin{array}{c}0\\1/2\end{array}\right](z)\right)^2.
\left(\theta\left[\begin{array}{c}0\\1/2\end{array}\right](0)\right)^2\nonumber\\
&&+
\left(\theta\left[\begin{array}{c}1/2\\0\end{array}\right](z)\right)^2.
\left(\theta\left[\begin{array}{c}1/2\\0\end{array}\right](0)\right)^2,\nonumber
\end{eqnarray}
et
\begin{eqnarray}
\left(\theta\left[\begin{array}{c}1/2\\1/2\end{array}\right](z)\right)^2.
\left(\theta\left[\begin{array}{c}0\\0\end{array}\right](0)\right)^2
&=&
\left(\theta\left[\begin{array}{c}0\\1/2\end{array}\right](z)\right)^2.
\left(\theta\left[\begin{array}{c}1/2\\0\end{array}\right](0)\right)^2\nonumber\\
&&-
\left(\theta\left[\begin{array}{c}1/2\\0\end{array}\right](z)\right)^2.
\left(\theta\left[\begin{array}{c}0\\1/2\end{array}\right](0)\right)^2,\nonumber
\end{eqnarray}
ainsi que l'identit\'{e} de Jacobi obtenue en posant $z=0$,
$$\left(\theta\left[\begin{array}{c}0\\0\end{array}\right](0)\right)^4
=\left(\theta\left[\begin{array}{c}0\\1/2\end{array}\right](0)\right)^4
+\left(\theta\left[\begin{array}{c}1/2\\0\end{array}\right](0)\right)^4.$$

Nous allons voir comment exprimer les fonctions m\'{e}romorphes
sur le tore $\mathbb{C}/\Lambda$, en termes de la fonction
th\^{e}ta. Plusieurs approches sont possibles :

\emph{\underline{Approche 1} :} Rappelons que toute fraction
rationnelle (donc fonction m\'{e}romorphe sur
$\mathbb{P}^1(\mathbb{C})$) peut s'\'{e}crire sous la forme
$$f(z)=\prod_{j=1}^m\frac{z-P_j}{z-Q_j}.$$
Par analogie, soient $P_1,...,P_m,Q_1,...,Q_m$ des points de la
surface de Riemann $X$ et $f(z)$ une fonction ayant des z\'{e}ros
aux points $P_1,...,P_m$ et des p\^{o}les aux points
$Q_1,...,Q_m$. On suppose que la condition (i) (ou ce qui est
\'{e}quivalent, la condition (ii)) du th\'{e}or\`{e}me
d'Abel\footnote{Soient $p_1,...,p_m, q_1,...q_m$ des points de
$X$. Alors les deux conditions suivantes sont \'{e}quivalentes :
(i) Il existe une fonction m\'{e}romorphe $f$ telle que : $(f)
=\sum_{j=1}^{m}q_{j}-\sum_{j=1}^{m}p_{j}$. (ii) Il existe un
chemin ferm\'{e} $\gamma $\ tel que : $\forall \omega \in
\Omega^1(X) , \quad\sum_{j=1}^{m}\int_{p_{j}}^{q_{j}}\omega
=\int_{\gamma }\omega$. Soit $\mathcal{D}=\sum_{j=1}^mn_jq_j\in
\mbox{Div }(X)$, $p\in X$, fix\'{e} et soit
$(\omega_1,...,\omega_g)$ une base de diff\'{e}rentielles
holomorphes sur $X$. L'application
$$\varphi : \mbox{Div }(X)\longrightarrow \mbox{Jac}(X),\quad
\mathcal{D}\longmapsto
\left(\sum_{j=1}^{m}n_j\int_{p}^{q_{j}}\omega_1,...,\sum_{j=1}^{m}n_j\int_{p}^{q_{j}}\omega_g\right),$$
est dite "application d'Abel-Jacobi". En particulier, si
$\mathcal{D}=\mathcal{D}_1-\mathcal{D}_2=\sum_{j=1}^mq_j-\sum_{j=1}^mp_j$,
alors la condition (i) signifie que $\mathcal{D}\in
\mbox{Div}^0(X)$ ou encore $\mathcal{D}_1$ est \'{e}quivalent
\`{a} $\mathcal{D}_2$. La condition (ii) peut s'\'{e}crire sous
une forme condens\'{e}e, $\forall\omega \in \Omega^1(X)$,
$\int_{\mathcal{D}_1}^{\mathcal{D}_2}\omega =\int_{\gamma}\omega$
ou encore sous la forme $\varphi(\mathcal{D})\equiv
\left(\sum_{j=1}^{m}\int_{p_{j}}^{q_{j}}\omega_1,...,\sum_{j=1}^{m}\int_{p_{j}}^{q_{j}}\omega_g\right)
\equiv 0 \mbox{ mod. } L$, avec $\varphi$ l'application
d\'{e}finie par $\varphi : \mbox{Div}^\circ (X)\longrightarrow
\mbox{Jac}(X)$.} est satisfaite. Comme $X$ est de genre $1$, alors
il existe une seule diff\'{e}rentielle holomorphe $\omega$ sur
$X$. Toujours d'apr\`{e}s le th\'{e}or\`{e}me d'Abel, l'existence
de la fonction $f(z)$ impose la condition
$\sum_{j=1}^mP_j=\sum_{j=1}^mQ_j$. Notons que pour $m=1$;
$P_1=Q_1$ et le seul cas valable est $f(z)=\mbox{constante}$. Dans
le cas o\`{u} $m\geq2$, alors la fonction $f(z)$ s'exprime en
fonction de $\theta$ \`{a} l'aide de la formule
$$f(z)=C\prod_{j=1}^m\frac{\theta\left(z-P_j-\frac{1}{2}(1+b)\right)}{\theta\left(z-Q_j-\frac{1}{2}(1+b)\right)},$$
o\`{u} $C$ est une constante. Notons que $f(z+1)=f(z)$. En outre,
d'apr\`{e}s la relation (1) et du fait que
$\sum_{j=1}^mP_j=\sum_{j=1}^mQ_j$, on a aussi $f(z+b)=f(z)$. Donc
$f$ est doublement p\'{e}riodique. La fonction $f$ est
m\'{e}romorphe avec des z\'{e}ros en $Q_j+\frac{1}{2}(1+b)$ et des
p\^{o}les en $P_j+\frac{1}{2}(1+b)$.

\emph{\underline{Approche 2} :} La fonction $\log\theta(z)$ peut
s'exprimer comme \'{e}tant la somme d'une fonction doublement
p\'{e}riodique de p\'{e}riodes $1$, $b$ et d'une fonction
lin\'{e}aire. Donc la fonction $\frac{d^2}{dz^2}\log\theta(z)$ est
doublement p\'{e}riodique et m\'{e}romorphe sur $X$, avec un
p\^{o}le double en $z=\frac{1}{2}(1+b)$. Cette fonction coincide
avec la fonction $\wp(z)$ de Weierstrass :
\begin{equation}\label{eqn:euler}
\wp(z)=-\frac{d^2}{dz^2}\log\theta(z)+C,
\end{equation}
o\`{u} $C$ est une constante choisie de telle mani\`{e}re que le
d\'{e}veloppement en s\'{e}rie de Laurent de $\wp(z)$ en $z=0$ n'a
pas de terme constant. On montre que la fonction $\theta(z)$
satisfait \`{a} une \'{e}quation diff\'{e}rentielle de
$3^{\mbox{ème}}$-ordre. Il suffit d'utiliser la relation (1) et
l'\'{e}quation diff\'{e}rentielle (6).

\emph{\underline{Approche 3} :} Rappelons que les fonctions
m\'{e}romorphes avec des p\^{o}les simples sur
$\mathbb{P}^1(\mathbb{C})$ peuvent s'\'{e}crire sous la forme
$$f(z)=\sum_j\frac{\lambda_j}{z-P_j}+C,$$
o\`{u} $\lambda_j\in\mathbb{C}$ et $C$ une constante. Par
analogie, on consid\`{e}re sur $X$ la fonction
$$f(z)=\sum_j\lambda_j\frac{d}{dz}\log\theta(z-P_j)+C,$$
o\`{u} $P_j\in X$, $\lambda_j\in\mathbb{C}$ tel que :
$\sum_j\lambda_j=0$ et $C$ est une constante. Cette fonction est
doublement p\'{e}riodique et m\'{e}romorphe avec des p\^{o}les
simples en $P_j+\frac{1}{2}(1+b)$ et de r\'{e}sidus $\lambda_j$ en
ces points.

Nous avons vu comment peuvent s'exprimer les fonctions
m\'{e}romorphes sur le tore $\mathbb{C}/\Lambda$ en termes de
fonction th\^{e}ta. Par ailleurs, pour $g=1$, on sait que :
$X\simeq\mathbb{C}/\Lambda\simeq\mbox{Jac}(X)$. Donc la
construction qui a \'{e}t\'{e} faite pr\'{e}c\'{e}demment sur le
tore $\mathbb{C}/\Lambda$ ou ce qui revient au m\^{e}me sur
$\mbox{Jac}(X)$ est aussi valable sur la surface de Riemann $X$.
Par exemple, prenons le cas d'une fonction ayant des p\^{o}les en
$P_1,...,P_m$ et des z\'{e}ros en $Q_1,...,Q_m$ sur la surface de
Riemann $X$. D'apr\`{e}s le th\'{e}or\`{e}me d'Abel, on a
$\sum_{j=1}^m\varphi(P_j)=\sum_{j=1}^m\varphi(Q_j)$, et on peut
selon la m\'{e}thode 1 d\'{e}crite ci-dessus, exprimer la fonction
$f(P)$ en termes de fonction th\^{e}ta à l'aide de la formule
$$f(P)=C\prod_{j=1}^m\frac{\theta\left(\varphi(P)-\varphi(Q_j)-\frac{1}{2}(1+b)\right)}
{\theta\left(\varphi(P)-\varphi(P_j)-\frac{1}{2}(1+b)\right)}.$$

Passons maintenant au cas o\`{u} la surface de Riemann $X$ est de
genre $g>1$. Rappelons que le probl\`{e}me d'inversion de Jacobi
[5], consiste \`{a} d\'{e}terminer $g$ points $P_1,...,P_g$ sur
$X$ tels que :
$$\sum_{k=1}^g\int_{P_0}^{P_k}\omega_j\equiv z_j\quad(\mbox{mod
}L),\quad j=1,...,g$$ o\`{u} $(z_1,...,z_g)\in\mbox{Jac}(X)$,
$(\omega_1,...,\omega_g)$ une base de diff\'{e}rentielles
holomorphes sur $X$, $P_0$ un point de base sur $X$ et $L$ un
r\'{e}seau engendr\'{e} par les vecteurs colonnes de la matrice
des p\'{e}riodes. Autrement dit, le probl\`{e}me consiste \`{a}
d\'{e}terminer le diviseur $\mathcal{D}=\sum_{j=1}^gP_j$ en termes
de $z=(z_1,...,z_g)\in\mbox{Jac}(X)$ tel que si $\varphi$ est
l'application d'Abel-Jacobi, alors l'\'{e}quation
$\varphi(\mathcal{D})=z$ soit satisfaite. Nous allons \'{e}tudier
le probl\`{e}me d'inversion de Jacobi \`{a} l'aide des fonctions
th\^{e}ta.

\begin{thm}
Si la fonction d\'{e}finie par
$\zeta(P)=\theta\left(\varphi(P)-C\right)$, $C\in\mathbb{C}^g$,
n'est pas identiquement nulle, alors elle admet $g$ z\'{e}ros
(compt\'{e}s avec leur ordre de multiplicit\'{e}) sur la
repr\'{e}sentation normale $X^*$ de $X$, que l'on d\'{e}signe par
le symbole $a_1b_1a_1^{-1}b_1^{-1}\ldots a_gb_ga_g^{-1}b_g^{-1}$,
o\`{u} $(a_{1},\ldots ,a_{g},b_{1},\ldots ,b_{g})$ est une base
symplectique du groupe d'homologie $H_{1}( X,\mathbb{Z})$. En
outre, si $P_1,...,P_g$ d\'{e}signent les z\'{e}ros de cette
fonction alors on a sur la vari\'{e}t\'{e} jacobienne
$\mbox{Jac}(X)$ la formule
$\displaystyle{\sum_{k=1}^g\varphi(P_k)\equiv C-\Delta}$, (mod.
p\'{e}riodes), o\`{u} $\Delta\in\mathbb{C}^g$ est le vecteur
des constantes de Riemann d\'{e}fini par
\begin{equation}\label{eqn:euler}
\Delta_j=\frac{1}{2}(1+b_{jj})-\sum_{k\neq
j}\left(\int_{a_k}\omega_k(P)\int_{P_0}^P\omega_j\right),\quad
j=1,...,g.
\end{equation}
\end{thm}
\emph{D\'{e}monstration}: Notons que $X^*$ est un polyg\^{o}ne à $4g$ c\^{o}t\'{e}s identifi\'{e}s deux
\`{a} deux. Si l'on parcourt le bord $\partial X^*$ de ce
polyg\^{o}ne, on constate que chaque c\^{o}t\'{e} est parcouru
deux fois, l'un dans le sens de son orientation et l'autre dans le
sens oppos\'{e}. On a donc $\displaystyle{\partial
X^*=\sum_{j=1}^g\left(a_j+b_j-a_j^{-1}-b_j^{-1}\right)}$. On
d\'{e}signe par $\zeta^-$ la valeur de la fonction $\zeta(P)$ sur
$a_j^{-1}$, $b_j^{-1}$ et par $\zeta^+$ la valeur de $\zeta(P)$
sur les segments $a_j$, $b_j$. On utilisera des notations
similaires $\varphi^+$, $\varphi^-$ pour l'application d'Abel
$\varphi(P)$. Le nombre de z\'{e}ros de la fonction $\zeta$ sur
$X^*$ est
\begin{equation}\label{eqn:euler}
\frac{1}{2\pi i}\int_{\partial X^*}d\log\zeta(P) =\frac{1}{2\pi
i}\sum_{k=1}^g\left(\int_{a_k}+\int_{b_k}\right)(d\log\zeta^+-d\log\zeta^-).
\end{equation}
Notons que : $\varphi_j^-(P)=\varphi_j^+(P)+b_{jk}$ si $P\in a_k$
et $\varphi_j^+(P)=\varphi_j^-(P)+\delta_{jk}$ si $P\in b_k$.
D'apr\`{e}s le th\'{e}or\`{e}me 4, on a
\begin{eqnarray}
d\log\varphi^-(P)&=&d\log\varphi^+(P)-2\pi i\omega_k  \mbox{
sur }a_k,\nonumber\\
d\log\varphi^+(P)&=&d\log\varphi^-(P) \mbox{ sur }b_k.\nonumber
\end{eqnarray}
Par cons\'{e}quent, (2) implique
$$\frac{1}{2\pi i}\int_{\partial X^*}d\log\zeta=\frac{1}{2\pi
i}\sum_{k=1}^g\int_{a_k}2\pi i\omega_k=g,$$ ce qui montre que la
fonction $\zeta(P)$ admet $g$ z\'{e}ros sur $X^*$. Pour prouver la
seconde partie du th\'{e}or\`{e}me, on consid\`{e}re
l'int\'{e}grale
$$I_j=\int_{\partial X^*}\varphi_j(P)d\log\zeta(P),\quad
j=1,...,g.$$ En d\'{e}signant par $P_1,...,P_g$ les z\'{e}ros de
la fonction $\zeta(P)$ et en tenant compte du th\'{e}or\`{e}me des
r\'{e}sidus, on a $I_j=\varphi_j(P_1)+\cdots+\varphi_j(P_g)$. En
raisonnant comme pr\'{e}c\'{e}demment, on obtient
\begin{eqnarray}
I_j&=&\frac{1}{2\pi
i}\sum_{k=1}^g\left(\int_{a_k}+\int_{b_k}\right)
\left(\varphi_j^+d\log\zeta^+-\varphi_j^-d\log\zeta^-\right),\nonumber\\
&=&\frac{1}{2\pi
i}\sum_{k=1}^g\int_{a_k}\left(\varphi_j^+d\log\zeta^+-(\varphi_j^++b_{jk})(d\log\zeta^+-2\pi
i\omega_k)\right)\nonumber\\
&&+\frac{1}{2\pi
i}\sum_{k=1}^g\int_{b_k}\left(\varphi_j^+d\log\zeta^+-(\varphi_j^+-\delta_{jk})d\log\zeta^+\right),\nonumber\\
&=&\sum_{k=1}^g\left(\int_{a_k}\varphi_j^+\omega_k-\frac{1}{2\pi
i}b_{jk}\int_{a_k}d\log\zeta^++b_{jk}\right)+\frac{1}{2\pi
i}\int_{b_k}d\log\zeta^+ .\nonumber
\end{eqnarray}
Notons que $\int_{a_k}d\log\zeta^+=2\pi in_k$, $n_k\in
\mathbb{Z}$. De m\^{e}me, en d\'{e}signant par $Q_j$ (resp.
$Q_j^*$) le d\'{e}but (resp. fin) du contour $b_j$, alors
\begin{eqnarray}
\int_{b_j}d\log\zeta^+&=&\log\zeta^+(Q_j^*)-\log\zeta^+(Q_j)+2\pi
im_j,\quad m_j\in\mathbb{Z},\nonumber\\
&=&\log \theta\left(\varphi(Q_j)+f_j-C\right)-\log
\theta\left(\varphi(Q_j)-C\right)+2\pi im_j,\nonumber\\
&=&-\pi ib+2\pi iC_j-2\pi i\varphi_j(Q_j)+2\pi im_j,\nonumber
\end{eqnarray}
o\`{u} $f_j= (b_{1j} ... b_{gj})^\top$, $j=1,...,g$, d\'{e}signent
les colonnes de la matrice $B$. D\`{e}s lors,
$$I_j=C_j-\frac{1}{2}b_{jj}-\varphi_j(Q_j)+\sum_{k=1}^g\int_{a_k}\varphi_j(P)\omega_k\quad(\mbox{mod.
p\'{e}riodes}).$$ Le d\'{e}but du contour $a_j$ sera
d\'{e}sign\'{e} par $R_j$ et sa fin coincide \'{e}videmment avec
le d\'{e}but $Q_j$ du contour $b_j$. On a
\begin{eqnarray}
I_j&=&C_j-\frac{1}{2}b_{jj}-\varphi_j(Q_j)+\int_{a_j}\varphi_j(P)\omega_j
+\sum_{\underset{k\neq j}{k=1}}^g\int_{a_k}\varphi_j(P)\omega_k,\nonumber\\
&=&C_j-\frac{1}{2}b_{jj}-\varphi_j(Q_j)+\frac{1}{2}\left(\varphi_j^2(Q_j)-\varphi_j^2(R_j)\right)
+\sum_{\underset{k\neq j}{k=1}}^g\int_{a_k}\varphi_j(P)\omega_k,\nonumber\\
&=&C_j-\frac{1}{2}b_{jj}-\varphi_j(R_j)-1+\frac{1}{2}\left((\varphi_j^2(R_j)+1)^2-\varphi_j^2(R_j)\right)
+\sum_{\underset{k\neq j}{k=1}}^g\int_{a_k}\varphi_j(P)\omega_k,\nonumber\\
&=&C_j-\frac{1}{2}(1+b_{jj})+\sum_{\underset{k\neq
j}{k=1}}^g\int_{a_k}\varphi_j(P)\omega_k,\nonumber
\end{eqnarray}
ce qui termine la preuve. $\square$

En g\'{e}n\'{e}ral, le vecteur $\Delta$ d\'{e}pend de $P_0$ sauf
dans le cas particulier $g=1$ o\`{u} $\Delta=\frac{1}{2}(1+b)$. On
montre que $2\Delta=-\varphi(K)$, o\`{u} $K$ est le diviseur
canonique. D\`{e}s lors, en choisissant adroitement le point
$P_0$, on peut exprimer $K$ de mani\`{e}re tout \`{a} fait simple.
Par exemple, consid\'{e}rons le cas o\`{u} $X$ est une courbe
hyperelliptique de genre $g$ d'\'{e}quation affine
$\displaystyle{w^2=\prod_{j=1}^{2g+2}(\xi-\xi_j)}$, o\`{u} tous
les $z_j$ sont distincts. Soit $(a_{1},\ldots ,a_{g},b_{1},\ldots
,b_{g})$ une base symplectique du groupe d'homologie $H_{1}(
X,\mathbb{Z})$) et soit $\sigma :X\longrightarrow X$,
$(w,\xi)\longmapsto (-w,\xi)$, l'involution hyperelliptique
(c.-\`{a}-d., qui consiste \`{a} \'{e}changer les deux feuillets
de la courbe $X$) avec $\sigma(a_j)=-a_j$ et $\sigma(b_j)=-b_j$.
Notons que
$$\int_{a_j}\omega_k=-\int_{\sigma(a_j)}\omega_k=-\int_{a_j}\sigma^*\omega_k.$$
Alors, en choisissant $P_0=\xi_1$, on obtient
\begin{eqnarray}
\Delta_j&=&\frac{1}{2}(1+b_{jj})+\sum_{k\neq
j}\int_{a_k}\omega_k\left(\int_{\xi_1}^{\xi_{2k+1}}\omega_j+\int_{\xi_{2k+1}}^P\omega_j\right),\quad
j=1,...,g,\nonumber\\
&=&\frac{1}{2}(1+b_{jj})+\sum_{k\neq
j}\int_{\xi_1}^{\xi_{2k+1}}\omega_j\int_{a_k}\omega_k\nonumber\\
&&+\sum_{k\neq
j}\int_{\xi_{2k+1}}^{\xi_{2k+2}}\left(\left(\int_{\xi_{2k+1}}^P\omega_j\right)
\omega_k(P)\left(\int_{\xi_{2k+1}}^{\sigma
P}\omega_j\right)\omega_k(\sigma P)\right).\nonumber
\end{eqnarray}
En tenant compte du fait que $\omega_k(\sigma P)=-\omega_k(P)$ et
modulo une combinaison lin\'{e}aire $n+Bm$ (un r\'{e}seau
engendr\'{e} par les vecteurs colonnes de la matrice des
p\'{e}riodes), on obtient ce cas la formule :
$\Delta_j=\displaystyle{\sum_{k=1}^gb_{jk}+\frac{j}{2}}$, $1\leq
j\leq g$. Les z\'{e}ros d'une fonction th\^{e}ta sur
$\mathbb{C}^g$ forment une sous-vari\'{e}t\'{e} de $\mbox{Jac}(X)$
de dimension $g-1$ appel\'{e}e diviseur th\^{e}ta que l'on note
$\Theta=\{z :\theta(z)=0\}$. Elle est invariante par un nombre
fini de translations et peut \^{e}tre singuli\`{e}re.
L'\'{e}quation (3) implique que $\Theta$ est bien d\'{e}finie sur
la vari\'{e}t\'{e} jacobienne $\mbox{Jac}(X)$. Comme
$\theta(-z)=\theta(z)$, on en d\'{e}duit que $\Theta$ est
sym\'{e}trique : $-\Theta=\Theta$.

\begin{thm}
(Riemann [4]). La fonction
$\zeta(P)=\theta\left(\varphi(P)-C\right)$, $C\in\mathbb{C}^g$,
est soit identiquement nulle, soit admet exactement $g$ z\'{e}ros
$Q_1,...,Q_g$ sur $X$ tels que :
$\displaystyle{\sum_{j=1}^g\varphi(Q_j)=C+\Delta}$, o\`{u}
$\Delta$ est d\'{e}fini par (8).
\end{thm}
Ce r\'{e}sultat signifie que lorsqu'on
plonge la surface de Riemann $X$ dans sa vari\'{e}t\'{e}
jacobienne $\mbox{Jac}(X)$ via l'application $\varphi$ d'Abel,
alors soit son image est enti\`{e}rement inclue dans le diviseur
th\^{e}ta, soit elle la rencontre en exactement $g$ points. En
fait si $\zeta(P)$ n'est pas identiquement nulle sur $X$, alors
ses z\'{e}ros coincident avec les points $P_1,...,P_g$ et
d\'{e}terminent la solution du probl\`{e}me inverse de Jacobi
$\varphi(\mathcal{D})=z$ pour le vecteur $z=C-\Delta$. Rappelons
que $\mathcal{D}\in\mbox{Div}(X)$ est un diviseur sp\'{e}cial si
et seulement si $\dim\mathcal{L}(\mathcal{D})\geq 1$ et
$\dim\mathcal{L}(K-\mathcal{D})\geq 1$ o\`{u} $K$ est un diviseur
canonique. Dans le cas o\`{u} $\mathcal{D}\geq0$, un diviseur est sp\'{e}cial
si et seulement si $\dim\Omega^1(\mathcal{D})\neq 0$. Notons aussi
que les diviseurs sp\'{e}ciaux de la forme
$\mathcal{D}=P_1+\cdots+P_N$, $N=\mbox{deg }\mathcal{D}\geq g$,
coincident avec les points critiques de l'application
d'Abel-Jacobi,
$$\mbox{Sym}^NX\longrightarrow\mbox{Jac}(X),\quad
\mathcal{D}\longmapsto\left(\int_0^\mathcal{D}\omega_1,...,\int_0^\mathcal{D}\omega_N\right),$$
ou ce qui revient au m\^{e}me
$\varphi(P_1,...,P_N)=\varphi(P_1)+\cdots+\varphi(P_N)$. Ces
points critiques sont les points $P_1,...,P_N$ o\`{u} le rang de
la diff\'{e}rentielle de cette application est inf\'{e}rieur \`{a}
$g$. D'apr\`{e}s le th\'{e}or\`{e}me 7 fondamental de Riemann, la
fonction $\zeta(P)=\theta\left(\varphi(P)-C\right)$, est
identiquement nulle si et seulement si
$C\equiv\varphi(Q_1)+\cdots\cdots+\varphi(Q_g)+\Delta$ o\`{u}
$Q_1+\cdots+Q_g$ est un diviseur sp\'{e}cial.

\begin{thm}
Soit $z=(z_1,...,z_g)\in\mathbb{C}^g$ un vecteur tel que la
fonction $\zeta(P)=\theta\left(\varphi(P)-z-\Delta\right)$, n'est
pas identiquement nulle sur $X$. Alors, la fonction $\zeta(P)$
admet exactement $g$ z\'{e}ros $P_1,...,P_g$ sur $X$ qui
d\'{e}terminent la solution du probl\`{e}me d'inversion de Jacobi
$\varphi(\mathcal{D})=z$, o\`{u} $\mathcal{D}=\displaystyle{\sum_{j=1}^gP_j}$.
Autrement dit, on a
\begin{equation}\label{eqn:euler}
\varphi_j(P_1)+\cdots+\varphi_j(P_g)=\sum_{k=1}^g\int_{P_0}^{P_k}\omega_j\equiv
z_j,\quad 1\leq j\leq g
\end{equation}
En outre, le diviseur $\mathcal{D}$ est non sp\'{e}cial et les
points $P_1,...,P_g$ sont uniquement d\'{e}termin\'{e}s \`{a}
partir du syst\`{e}me (10).
\end{thm}
\emph{D\'{e}monstration}: La premi\`{e}re assertion r\'{e}sulte du
th\'{e}or\`{e}me 1. En outre, le diviseur
$\mathcal{D}=\displaystyle{\sum_{j=1}^gP_j}$ n'est pas sp\'{e}cial
car sinon la fonction $\zeta(P)$ serait identiquement nulle
d'apr\`{e}s ce qui pr\'{e}c\`{e}de, ce qui est absurde. Quand au
dernier point, supposons que le syst\`{e}me (10) admet une autre
solution $Q_1,...,Q_g$. On aura sur la vari\'{e}t\'{e} jacobienne
$\mbox{Jac}(X)$,
$\displaystyle{\sum_{j=1}^g\varphi(P_j)\equiv\sum_{j=1}^g\varphi(Q_j)}$,
(mod. $L$), o\`{u} $L$ est le r\'{e}seau engendr\'{e} par la
matrice des p\'{e}riodes. D'apr\`{e}s le th\'{e}or\`{e}me d'Abel,
cela signifie qu'il existe une fonction m\'{e}romorphe sur $X$
ayant des z\'{e}ros en $Q_1,...,Q_g$ et des p\^{o}les en
$P_1,...,P_g$. Or on vient de montrer que le diviseur est non
sp\'{e}cial, donc une telle fonction doit \^{e}tre une constante,
ce qui entraine que $P_j=Q_j$, $j=1,...,g$. $\square$

Par exemple, si $\mathcal{D}=\displaystyle{\sum_{j=1}^gP_j}$ est
un diviseur non sp\'{e}cial sur une surface de Riemann $X$ de
genre $g$, alors la fonction $\theta(\varphi
(P)-\varphi(\mathcal{D})-\Delta)$, admet exactement $g$ z\'{e}ros
sur $X$ aux points $P=P_1,...,P_g$.

On a la caract\'{e}risation suivante du diviseur th\^{e}ta :

\begin{thm}
On a $\theta(C)=0$, si et seulement s'il existe
$P_1,...,P_{g-1}\in X$ avec point de base $P_0$, tels que :
$$
C\equiv\varphi(P_1)+\cdots+\varphi(P_{g-1})+\Delta=\sum_{j=1}^{g-1}\int_{P_0}^{P_j}\omega+\Delta
$$
\end{thm}
\emph{D\'{e}monstration}: Reprenons la fonction
$\zeta(P)=\theta(\varphi(P)-C)$ et supposons d'abord qu'elle est
non nulle sur $X$. D'apr\`{e}s le th\'{e}or\`{e}me 6, cette
fonction admet $g$ z\'{e}ros $P_1,...,P_g$ sur $X$ et
\begin{equation}\label{eqn:euler}
C\equiv\varphi(P_1)+\cdots+\varphi(P_{g})+\Delta.
\end{equation}
L'ensemble de ces z\'{e}ros \'{e}tant unique et comme par
hypoth\`{e}se $\theta(C)=0$, alors $P_g=P_0$. D\`{e}s lors
$\varphi(P_g)=\varphi(P_0)=0$ et d'apr\`{e}s (11), on a
$$C\equiv\varphi(P_1)+\cdots+\varphi(P_{g-1})+\Delta.$$
Passons maintenant au cas o\`{u} la fonction $\zeta(P)$ n'est pas
identiquement nulle sur $X$. D'apr\`{e}s le th\'{e}or\`{e}me 6, on
a
\begin{equation}\label{eqn:euler}
C\equiv\varphi(Q_1)+\cdots+\varphi(Q_{g})+\Delta,
\end{equation}
o\`{u} $Q_1+\cdots+Q_g$ est un diviseur sp\'{e}cial. Ce dernier
implique l'existence sur $X$ d'une fonction non constante $\zeta$
m\'{e}romorphe ayant des p\^{o}les en $Q_1,...,Q_g$ avec
$\zeta(P_0)=0$. D\`{e}s lors,
$\varphi(P_1+\cdots+P_{g-1}+P_0)\equiv\varphi(Q_1+\cdots+Q_g)$, en
vertu du th\'{e}or\`{e}me d'Abel o\`{u} $P_1+\cdots+P_{g-1}+P_0$
est le diviseur des z\'{e}ros de $\zeta$. Il suffit d\`{e}s lors
de remplacer dans (12), $\varphi(Q_1+\cdots+Q_g)$ par
$\varphi(P_1+\cdots+P_{g-1}+P_0)$ tout en tenant compte du fait
que $\varphi(P_0)=0$. $\square$

\begin{thm}
Soient $\mathcal{D}$ un diviseur non sp\'{e}cial de degr\'{e} $g$,
$\mathcal{D}'$ un diviseur positif de degr\'{e} $n$,
$(\omega_1,...,\omega_g)$ une base de diff\'{e}rentielles
holomorphes sur $X$,
$\varphi(P)=\left(\int_{P_0}^P\omega_1,...,\int_{P_0}^P\omega_g\right)$
l'application d'Abel de point de base $P_0$, $\eta$ une
diff\'{e}rentielle normalis\'{e}e de $3^{\mbox{\`{e}me}}$
esp\`{e}ce\footnote{Une diff\'{e}rentielle (ab\'{e}lienne)
m\'{e}romorphe sur $X$ est dite de $3^{\grave{e}me}$esp\`{e}ce si
elle n'a que des p\^{o}les simples (et son r\'{e}sidu est non nul
en au moins un point de $X$).} sur $X$ ayant des p\^{o}les sur
$\mathcal{D}'$ et de r\'{e}sidus $-1$, $U=(U_1,...,U_g)$ le
vecteur des p\'{e}riodes avec $U_k=\int_{b_k}\eta$ et enfin
$\Delta$ le vecteur d\'{e}fini \`{a} l'aide des constantes de
Riemann par (8). Si $\psi$ est une fonction m\'{e}romorphe sur $X$
ayant $g+n$ p\^{o}les sur $\mathcal{D}+\mathcal{D}'$, alors cette
fonction s'exprime en termes de fonction th\^{e}ta \`{a} l'aide de
la formule
$$\psi(P)=A\frac{\theta\left(\varphi(P)-\varphi(\mathcal{D})+U-\Delta\right)}
{\theta\left(\varphi(P)-\varphi(\mathcal{D})-\Delta\right)}e^{\int_{P_0}^P\eta},\quad A=\mbox{constante}.$$
\end{thm}
\emph{D\'{e}monstration}: Il est \`{a} noter que le contour
d'int\'{e}gration dans les int\'{e}grales $\int_{P_0}^P\eta$ et
$\int_{P_0}^P\omega_j$, $j=1,...,g$, est le m\^{e}me. La fonction
$\psi(P)$ admet des p\^{o}les uniquement sur
$\mathcal{D}+\mathcal{D}'$. Montrons que cette fonction est bien
d\'{e}finie sur $X$; c.-\`{a}-d., elle ne d\'{e}pend pas du chemin
d'int\'{e}gration. Autrement dit qu'elle ne change pas lorsque $P$
parcourt un cycle quelconque
$\gamma=\displaystyle{\sum_{k=1}^g(n_ka_k+m_kb_k)\in
H_1(X,\mathbb{Z})}$, $(n_k, m_k\in\mathbb{Z})$. Les expressions
$\int_{P_0}^P\eta$ et
$\varphi(P)=\left(\int_{P_0}^P\omega_1,...,\int_{P_0}^P\omega_g\right)$
se transforment respectivement comme suit :
$$\int_{P_0}^P\eta+\sum_{k=1}^gm_k\int_{b_k}\eta=\int_{P_0}^P\eta+2
i\langle m,U\rangle,\quad m=(m_1,...,m_g)\in\mathbb{Z}^g,$$ et
$\varphi(P)\longmapsto \varphi(P)+n+Bm$,
$n=(n_1,...,n_g)\in\mathbb{Z}^g$. Par ailleurs, en utilisant la
formule (4), on obtient
$$
\frac{\theta\left(\varphi(P)-\varphi(\mathcal{D})+U-\Delta\right)}
{\theta\left(\varphi(P)-\varphi(\mathcal{D})-\Delta\right)}=
\frac{e^{-\pi i\langle Bm,m\rangle-2\pi i\langle
m,\varphi(P)-\varphi(\mathcal{D})+U-\Delta\rangle}} {e^{-\pi
i\langle Bm,m\rangle-2\pi i\langle
m,\varphi(P)-\varphi(\mathcal{D})-\Delta\rangle}}=e^{-2\pi
i\langle m,U\rangle},
$$
et le r\'{e}sultat d\'{e}coule de la transformation ci-dessus.
$\square$

Sur la surface de Riemann $X$ de genre $g$, des fonctions
singuli\`{e}res poss\`{e}dant $g$ p\^{o}les et des
singularit\'{e}s essentielles, jouent un r\^{o}le crucial lors de
l'\'{e}tude des syst\`{e}mes int\'{e}grables, notamment
l'\'{e}quation de Korteweg-de Vries (K-dV), $\frac{\partial
u}{\partial t}-6u\frac{\partial u}{\partial x}+\frac{\partial^3
u}{\partial x^3}=0$, l'\'{e}quation Kadomtsev-Petviashvili (KP),
$\frac{\partial^2 u}{\partial y^2}-\frac{\partial}{\partial
x}\left(4\frac{\partial u}{\partial t}-12u\frac{\partial
u}{\partial x}-\frac{\partial^3 u}{\partial x^3}\right)=0$,
l'\'{e}quation de Schr\"{o}dinger nonlin\'{e}aire $i\frac{\partial
\psi}{\partial t}+\frac{\partial^2\psi}{\partial
x^2}+\psi|^2\psi=0$, l'\'{e}quation de Boussinesq
$\frac{\partial^2u}{\partial t^2}-\frac{\partial^2u}{\partial
x^2}+\frac{\partial^4u}{\partial x^4}
+\frac{\partial^2u^2}{\partial x^2}=0$, l'\'{e}quation de
Camassa-Holm $\frac{\partial u}{\partial t}-\frac{\partial^3
u}{\partial t\partial x^2}+3u\frac{\partial u}{\partial
x}=2\frac{\partial u}{\partial x}\frac{\partial^2 u}{\partial
x^2}+u\frac{\partial^3 u}{\partial x^3}$, dont les solutions
exactes sont des solitons [14], c.-\`{a}-d., des ondes de formes
d\'{e}finies progressant \`{a} des vitesses diff\'{e}rentes. Nous
allons voir par analogie au th\'{e}or\`{e}me pr\'{e}c\'{e}dent,
comment exprimer ces fonctions (connues sous le nom de fonctions
de Baker-Akhiezer) en termes de fonctions th\^{e}ta et en m\^{e}me
temps prouver leur existence. Soient $Q_1,...,Q_n$ des points sur
une surface de Riemann $X$ de genre $g$ et $z_j$ des
param\`{e}tres locaux tels que : $z_j(Q_j)=\infty$. On associe
\`{a} chaque point $Q_j$ un polyn\^{o}me arbitraire not\'{e}
$q_j(z_j)$. Soient $\mathcal{D}=P_1+\cdots+P_g$ un diviseur
positive sur $X$ et $\psi(P)$ une fonction (dite fonction de
Baker-Akhiezer\index{fonction de Baker-Akhiezer}) satisfaisant aux
conditions suivantes : $(i)$ $\psi(P)$ est m\'{e}romorphe sur
$X\setminus\{Q_1,...,Q_n\}$ et admet des p\^{o}les uniquement aux
points $P_1,...,P_n$ du diviseur $\mathcal{D}$. $(ii)$ La fonction
$\psi(P)e^{-q_j(z_j(P))}$ est analytique au voisinage de $Q_j$,
$j=1,...,n$. On peut remplacer la condition $(ii)$ par celle-ci :
la fonction $\psi$ admet aux points $Q_j$, $j=1,...,n$, une
singularit\'{e} essentielle de la forme $\psi(P)\sim
ce^{q_j(z_j(P))}$ o\`{u} $c$ est une constante. Ces fonctions
$\psi(P)$ forment un espace vectoriel que l'on note $L\equiv
L(\mathcal{D};Q_1,...,Q_n,q_1,...,q_n)$.

\begin{thm}
Soit $\mathcal{D}=P_1+\cdots+P_g$ un diviseur non sp\'{e}cial de
degr\'{e} $g$. Alors l'espace $L$ est de dimension $1$ et sa base
est d\'{e}crite \`{a} l'aide de
\begin{equation}\label{eqn:euler}
\psi_1(P)=\frac{\theta\left(\varphi(P)-\varphi(\mathcal{D})+V-\Delta\right)}
{\theta\left(\varphi(P)-\varphi(\mathcal{D})-\Delta\right)}e^{\int_{P_0}^P\eta},
\end{equation}
o\`{u} $\eta$ est une diff\'{e}rentielle normalis\'{e}e de
$2^{\grave{e}me}$esp\`{e}ce\footnote{Une diff\'{e}rentielle
(ab\'{e}lienne) m\'{e}romorphe sur $X$ est dite de
$2^{\grave{e}me}$esp\`{e}ce si elle a des p\^{o}les et si son
r\'{e}sidu est nul en chaque point de $X$.} ayant des p\^{o}les
aux points $Q_1,...,Q_n$, les parties principales coincident avec
les polyn\^{o}mes $q_j(z_j)$, o\`{u} $j=1,...,n$,
$V=(V_1,...,V_g)$ avec $V_k=\int_{b_k}\eta$, $k=1,...,g$,
$$\varphi(P)=\left(\int_{P_0}^P\omega_1,...,\int_{P_0}^P\omega_g\right),$$
l'application d'Abel de point de base $P_0$, $\Delta$ est le
vecteur d\'{e}fini \`{a} l'aide des constantes de Riemann par (9).
Le contour d'int\'{e}gration dans les int\'{e}grales
$\int_{P_0}^P\eta$ et $\int_{P_0}^P\omega_j$, $j=1,...,g$ est le
m\^{e}me.
\end{thm}
\emph{D\'{e}monstration}: La fonction $\psi_1(P)$ poss\`{e}de des
p\^{o}les sur le diviseur $\mathcal{D}$ et des singularit\'{e}s
essentielles aux points $Q_1,...,Q_n$. La fonction $\psi_1(P)$ est
bien d\'{e}finie; elle ne d\'{e}pend pas du chemin
d'int\'{e}gration. En utilisant les notations et un raisonnement
similaire \`{a} ceux du th\'{e}or\`{e}me 10, on obtient le rapport
$$\frac{\theta\left(\varphi(P)-\varphi(\mathcal{D})+V-\Delta\right)}
{\theta\left(\varphi(P)-\varphi(\mathcal{D})-\Delta\right)}=e^{-2\pi
i\langle m,V\rangle},$$ et le r\'{e}sultat d\'{e}coule de la
transformation utilis\'{e}e dans la preuve du th\'{e}or\`{e}me
pr\'{e}c\'{e}dent. Par ailleurs, d'apr\`{e}s le th\'{e}or\`{e}me
de Riemann-Rock [5], la dimension de l'espace $L$ est \'{e}gale
\`{a} $\mbox{deg }\mathcal{D}-g+1$. Comme $\mbox{deg
}\mathcal{D}=g$, alors la dimension de l'espace en question est
\'{e}gal à $1$, ce qui prouve l'unicit\'{e} de la fonction
$\psi_1$ \`{a} une constante multiplicative pr\`{e}s. Soit
$\psi\in L$ une fonction quelconque. D\`{e}s lors, le quotient
$\displaystyle{\frac{\psi}{\psi_1}}$ est une fonction
m\'{e}romorphe avec $g(=\mbox{deg }\mathcal{D})$ p\^{o}les. Le
diviseur des p\^{o}les de $\displaystyle{\frac{\psi}{\psi_1}}$
coincide avec le diviseur $\mathcal{D}'=P_1'+\cdots+P_g'$ des
z\'{e}ros de $\psi_1(P)$ et on doit avoir
$\varphi(\mathcal{D}')-\varphi(\mathcal{D})=V$. En choisissant les
polyn\^{o}mes $q_j$ avec des coefficients suffisamment petits ou
ce qui revient au m\^{e}me, les vecteurs de $V$ suffisamment
petits, alors la fonction th\^{e}ta qui se trouve dans le
num\'{e}rateur de l'expression ci-dessus n'est pas identiquement
nulle. Par cons\'{e}quent, son diviseur des p\^{o}les
$\mathcal{D}'$ n'est pas sp\'{e}cial et donc
$\displaystyle{\frac{\psi}{\psi_1}}$ est une constante. $\square$

\section{Exemples}

Il est bien connu que les solutions de nombreux syst\`{e}mes
int\'{e}grables sont donn\'{e}es en termes de fonctions th\^{e}ta
associ\'{e}es \`{a} des surfaces de Riemann compactes. Nous
verrons ci-dessous de telles solutions pour certains
probl\`{e}mes.

Comme premier exemple, on consid\`{e}re le mouvement d'un solide
dans un fluide parfait d\'{e}crit \`{a} l'aide des \'{e}quations
de Kirchhoff [3] :
\begin{equation}\label{eqn:euler}
\dot{p}=p\wedge \frac{\partial H}{\partial l},\qquad
\dot{l}=p\wedge \frac{\partial H}{\partial p}+l\wedge
\frac{\partial H}{\partial l},
\end{equation}
o\`{u} $p=(p_{1},p_{2},p_{3})\in \mathbb{R}^{3},$\
$l=(l_{1},l_{2},l_{3})\in \mathbb{R}^{3}$ et $H$\ l'hamiltonien.
Le syst\`{e}me (14) poss\`{e}de les trois int\'{e}grales
premi\`{e}res suivantes :
\begin{equation}\label{eqn:euler}
H_{1}=H,\qquad H_{2}=p_{1}^{2}+p_{2}^{2}+p_{3}^{2},\qquad
H_{3}=p_{1}l_{1}+p_{2}l_{2}+p_{3}l_{3}.
\end{equation}
On distingue deux cas int\'{e}grables : cas de Clebsch et cas de
Lyapunov-Steklov. Dans le cas de Clebsch, on a
$\displaystyle{H=\frac{1}{2}\sum_{k=1}^{3}\left(
a_{k}p_{k}^{2}+b_{k}l_{k}^{2}\right)}$, avec la condition
$(a_{2}-a_{3})b_{1}^{-1}+(a_{3}-a_{1})b_{2}^{-1}+(a_{1}
-a_{2})b_{3}^{-1}=0$. Le syst\`{e}me ci-dessus s'\'{e}crit sous la
forme d'un champ de vecteurs hamiltonien. Une quatri\`{e}me
int\'{e}grale premi\`{e}re est fournie par
\begin{equation}\label{eqn:euler}
H_{4}=\frac{1}{2}\sum_{k=1}^{3}\left( b_{k}p_{k}^{2}+\varrho
l_{k}^{2}\right) ,
\end{equation}
o\`{u} $\varrho$ est une constante satisfaisant \`{a}
$$\varrho=b_{1}(b_{2}-b_{3})(a_{2}-a_{3})^{-1}
=b_{2}(b_{3}-b_{1})(a_{3}-a_{1})^{-1}=b_{3}(b_{1}-b_{2})(a_{1}-a_{2})^{-1}.$$
La m\'{e}thode de r\'{e}solution obtenue par K\"{o}tter [9] est
ext\^{e}mement compliqu\'{e}e et repose sur un choix astucieux de
deux variables $s_{1}$ et $s_{2}$. En  utilisant la substitution
$b_{k}\longrightarrow \varrho b_{k}$, $1\leq k\leq 3$, et une
combinaison lin\'{e}aire appropri\'{e}e de $H_{1}$\ et $H_{2},$ on
peut r\'{e}ecrire les \'{e}quations précédentes sous la forme
$$
p_{1}^{2}+p_{2}^{2}+p_{3}^{2}=A,\qquad\qquad
b_{1}p_{1}^{2}+b_{2}p_{2}^{2}+b_{3}p_{3}^{2}+l_{1}^{2}+l_{2}^{2}+l_{3}^{2}=B,$$
$$b_{1}l_{1}^{2}+b_{2}l_{2}^{2}+b_{3}l_{3}^{2}-b_{2}b_{3}p_{1}^{2}-b_{1}b_{3}p_{2}^{2}-b_{1}b_{2}p_{3}^{2}=C,\qquad
p_{1}l_{1}+p_{2}l_{2}+p_{3}l_{3}=D,$$ o\`{u} $A$, $B$, $C$ et $D$
sont des constantes. Introduisons des coordonn\'{e}es $\varphi
_{k}$, $\psi _{k}$, $1\leq k\leq 3$, en posant
$\varphi_{k}=p_{k}T_{+1}+l_{k}S_{+1}$ et
$\psi_{k}=p_{k}T_{-1}+l_{k}S_{-1}$, o\`{u}
$$T_{\pm1}=\frac{\sqrt{\prod_{j=1}^{3}(z_{1}-b_{j})}}
{\sqrt{z_{1}-b_{k}} \sqrt{\frac{\partial R}{\partial z_{1}}}}
+i\frac{\sqrt{\prod_{j=1}^{3}
\left(z_{2}-b_{j}\right)}}{\sqrt{z_{2}-b_{k}}\sqrt{\frac{\partial
R}{\partial z_{2}}}}$$
$$
S_{\pm1}=\frac{\sqrt{z_{1}-b_{k}}}{\sqrt{\frac{\partial
R}{\partial z_{1}}}}
+i\frac{\sqrt{z_{2}-b_{k}}}{\sqrt{\frac{\partial R}{\partial
z_{2}}}},\qquad R\left( z\right) =\prod_{i=1}^{4}\left(
z-z_{i}\right),$$ et $z_{1},z_{2},z_{3},z_{4}$\ sont les racines
de l'\'{e}quation
$$A^{2}\left(z^{2}-z\sum_{k=1}^{3}b_{k}\right)+Bz-C+2D\sqrt{\prod_{k=1}^{3}\left(z-b_{k}\right)
}=0.$$ Soient $s_{1}$\ et $s_{2}$\ les racines de l'\'{e}quation
$$\psi _{1}^{2}\left(\nu _{1}^{2}-s\right)^{-1}+\psi _{2}^{2}\left(\nu _{2}^{2}-s\right)^{-1}+\psi_{3}^{2}\left(\nu
_{3}^{2}-s\right)^{-1}=0,$$ o\`{u}
$$\nu _{k}=\left(\frac{\sqrt{z_{3}-b_{k}}}{\sqrt{\frac{\partial R}{\partial z_{3}}}}
+i\frac{\sqrt{z_{4}-b_{k}}}{\sqrt{\frac{\partial R}{\partial
z_{4}}}}\right)\left(\frac{\sqrt{z_{1}-b_{k}}}
{\sqrt{\frac{\partial R}{\partial
z_{1}}}}+i\frac{\sqrt{z_{2}-b_{k}}}{\sqrt{\frac{\partial
R}{\partial z_{2}}}}\right)^{-1}, \quad 1\leq k\leq 3.$$ On peut
exprimer les variables $p_{1},p_{2},p_{3},l_{1},l_{2},l_{3}$ en
terme de $s_{1}$ et $s_{2}$ (voir [9]). Apr\`{e}s quelques
manipulations alg\'{e}briques, on obtient
$$
\dot{s}_{1}=\frac{\left( as_{1}+b\right) \sqrt{P_{5}\left(
s_{1}\right) }}{s_{2}-s_{1}},\qquad \dot{s}_{2}=\frac{\left(
as_{2}+b\right) \sqrt{P_{5}\left( s_{2}\right) }}{s_{1}-s_{2}},
$$
o\`{u} $a,b$ sont des constantes et $P_{5}\left( s\right)$ est un
polyn\^{o}me de degr\'{e} cinq ayant la forme suivante :
$P_{5}(s)=s(s-\nu_{1}^{2})(s-\nu_{2}^{2})(s-\nu_{3}^{2})(s-\nu_{1}^{2}\nu_{2}^{2}\nu_{3}^{2})$.
Par cons\'{e}quent, l'int\'{e}gration s'effectue au moyen de
fonctions hyperelliptiques de genre $2$ et les solutions peuvent
s'exprimer en termes de fonctions th\^{e}ta. Le probl\`{e}me de ce
mouvement est un cas limite du flot g\'{e}od\'{e}sique sur
$SO(4)$. Rappelons que pour un syst\`{e}me alg\'{e}briquement
compl\`{e}tement int\'{e}grable [1], on demande que les invariants
du syst\`{e}me diff\'{e}rentiel soient polynomiaux (dans des
coordonn\'{e}s ad\'{e}quates) et que de plus les vari\'{e}t\'{e}s
complexes obtenues en \'{e}galant ces invariants polynomiaux \`{a}
des constantes g\'{e}n\'{e}riques forment la partie affine d'un
tore complexe alg\'{e}brique (vari\'{e}t\'{e} ab\'{e}lienne) de
telle fa\c{c}on que les flots complexes engendr\'{e}s par les
invariants soient lin\'{e}aires sur ces tores complexes. Les
solutions m\'{e}romorphes d\'{e}pendant d'un nombre suffisant de
param\`{e}tres libres jouent un r\^{o}le crucial dans l'\'{e}tude
de ces systèmes. On montre [1, 6] que le syst\`{e}me
diff\'{e}rentiel en question est alg\'{e}briquement
compl\`{e}tement\ int\'{e}grable et le flot correspondant
\'{e}volue sur une surface ab\'{e}lienne $\widetilde{M_{c}}\cong
\mathbb{C}^{2}/L_{\Omega }$ o\`{u} le r\'{e}seau $L_{\Omega }$ est
engendr\'{e} par la matrice des p\'{e}riodes
$$\Omega=\left(\begin{array}{llll}
2 & 0 & a & c \\
0 & 4 & c & b
\end{array}
\right),\text{ Im}\left(
\begin{array}{ll}
a & c \\
c & b
\end{array}
\right)>0,\quad \left( a,b,c\in \mathbb{C}\right).$$ La surface
affine $M_{c}$ d\'{e}finie en \'{e}galant les invariants du
syst\`{e}me \`{a} des constantes g\'{e}n\'{e}riques, se
compl\`{e}te en $\widetilde{M_{c}}$ par l'adjonction d'une courbe
lisse $\mathcal{D}$ de genre $9$, laquelle est un rev\^{e}tement
ramifi\'{e} le long d'une courbe elliptique $\mathcal{E}$.
L'application $\widetilde{M_{c}}\longrightarrow\mathbb{CP}^7$,
$(t_1,t_2)\longmapsto[1,X_1(t_1,t_2),...,X_7(t_1,t_2)]$, est un
plongement de $\widetilde{M_{c}}$ dans l'espace projectif
$\mathbb{CP}^7$ où $(1,X_1,...,X_7)$ forme une base de l'espace
vectoriel $\mathcal{L}(\mathcal{D})$ des fonctions m\'{e}romorphes
ayant au plus un p\^{o}le simple sur $\mathcal{D}$ (les fonctions
$X_1,...,X_7$ s'exrime de mani\`{e}re simple en fonction de
$x_1,...,x_6$). Les solutions du syst\`{e}me diff\'{e}rentiel en
question sont donn\'{e}es en termes de fonctions th\^{e}ta par
$$X_k(t)=\frac{\theta_k[(t_1^0,t_2^o)+t(n_1,n_2)]}{\theta_0[(t_1^0,t_2^o)+t(n_1,n_2)]},\quad
k=1,...,7$$ où $(\theta_0,...,\theta_7)$ forme une base de
l'espace vectoriel des fonctions th\^{e}ta associ\'{e}es \^{e}
$\mathcal{D}$ (les deux fonctions th\^{e}ta $\theta_0$, $\theta_7$
sont impaires tandis que les six fonctions th\^{e}ta
$\theta_1$,..., $\theta_6$ sont paires. Pour le cas de
Lyapunov-Steklov, on a
$$H_{1}=H=\frac{1}{2}\sum_{k=1}^{3}\left( a_{k}p_{k}^{2}+b_{k}l_{k}^{2}\right)
+\sum_{k=1}^{3}c_{k}p_{k}l_{k},$$ $a_{1}=A^{2}b_{1}\left(
b_{2}-b_{3}\right) ^{2}+B$, $a_{2}=A^{2}b_{2}\left(
b_{3}-b_{1}\right) ^{2}+B$, $a_{3}=A^{2}b_{3}\left(
b_{1}-b_{2}\right) ^{2}+B$, $c_{1}=Ab_{2}b_{3}+C$,
$c_{2}=Ab_{1}b_{3}+C$, $c_{3}=Ab_{1}b_{2}+C$, o\`{u} $A$, $B$ et
$C$ sont des constantes. Une quatri\`{e}me int\'{e}grale
premi\`{e}re est fournie par
$$H_{4}=\frac{1}{2}\sum_{k=1}^{3}\left( d_{k}p_{k}^{2}+l_{k}^{2}\right)
-A\sum_{k=1}^{3}b_{k}p_{k}l_{k},$$ o\`{u} $ d_{1}=A^{2}\left(
b_{2}-b_{3}\right) ^{2}$, $d_{2}=A^{2}\left( b_{3}-b_{1}\right)
^{2}$, $d_{3}=A^{2}\left( b_{1}-b_{2}\right) ^{2}$. Un calcul long
[12] et d\'{e}licat montre que dans ce cas aussi,
l'int\'{e}gration s'effectue \`{a} l'aide de fonctions
hyperelliptiques de genre deux et les solutions peuvent s'exprimer
en termes de fonctions th\^{e}ta.

Un autre exemple concerne l'\'{e}quation de Landau-Lifshitz [2,
11] :
$$\frac{\partial S}{\partial t}=S\times\frac{\partial^2 S}{\partial x^2}+S\times JS,$$
o\`{u} $S=(S_1,S_2,S_3)$, $S_1^2+s_2^2+S_3^2=1$ et
$J=\mbox{diag}(J_1,J_2,J_3)$. Cette \'{e}quation d\'{e}crit les
effets d'un champ magn\'{e}tique sur les mat\'{e}riaux
ferromagn\'{e}tiques. Les solutions r\'{e}elles (avec
l'anisotropie magn\'{e}tique de type axe d'aimantation facile)
sont donn\'{e}es par
\begin{eqnarray}
S_1&=&\frac{\theta(\omega+d+m)\theta(\omega+d+m+r)-\theta(\omega+d)\theta(\omega+d+r)}{\theta(\omega+d)\theta(\omega+d+m+r)
-\theta(\omega+d+r)\theta(\omega+d+m)},\nonumber\\
S_2&=&-i\frac{\theta(\omega+d+m)\theta(\omega+d+m+r)+\theta(\omega+d)\theta(\omega+d+r)}{\theta(\omega+d)\theta(\omega+d+m+r)
-\theta(\omega+d+r)\theta(\omega+d+m)},\nonumber\\
S_3&=&\frac{\theta(\omega+d)\theta(\omega+d+m+r)+\theta(\omega+d+r)\theta(\omega+d+m)}{\theta(\omega+d)\theta(\omega+d+m+r)
-\theta(\omega+d+r)\theta(\omega+d+m)}.\nonumber
\end{eqnarray}
Ici la fonction theta est li\'{e}e \`{a} une courbe
hyperelliptique de genre $g$, le vecteur $d\in\mathbb{C}^g$ est
tel que : $\mbox{Im }d=-\frac{1}{2}\mbox{Im }r$,
$m=(m_1,...,m_g)$, $\omega=\frac{1}{2\pi}(Ux+Vt)$ et
$r=\int_{\bullet^+}^{\bullet^-}du$; le chemin de l'int\'{e}gration
doit \^ {e}tre fix\'{e}e.

Nous citons encore un autre exemple qui concerne l'\'{e}quation de
sine-Gordon [2] :
$$\frac{\partial^2 \varphi}{\partial x^2}-\frac{\partial^2 \varphi}{\partial t^2}=\sin\varphi.$$
C'est une \'{e}quation d'onde non-lin\'{e}aire aux applications
multiples en physique. Sa solution peut s'\'{e}crire sous la forme
$$\varphi(x,t)=2i\ln\frac{\theta\left[\begin{array}{c}\alpha\\\beta\end{array}\right](Ux+Vt+W|B)}{\theta\left[\begin{array}{c}0\\0\end{array}\right](Ux+Vt+W|B)}
+C+2\pi m,$$ o\`{u} $U,V,W\in\mathbb{C}^g, C\in\mathbb{R},
m\in\mathbb{Z}$.

Par ailleurs l'\'{e}tude des fonctions th\^{e}ta d'une surface de
Riemann du genre $g$ peut \^{e}tre faite \`{a} partir du point de
vue de la fonction tau d'une hi\'{e}rarchie d'\'{e}quations de
soliton [13]. Les fonctions tau sont des fonctions sp\'{e}cifiques
du temps, construites \`{a} partir de sections d'un fibr\'{e}
d\'{e}terminant sur une vari\'{e}t\'{e} grassmannienne de
dimension infinie et g\'{e}n\'{e}ralisent les fonctions th\^{e}ta
de Riemann.

\end{document}